\documentclass[a4paper,12pt]{article}
\usepackage{amssymb,amsmath,amsthm}
\usepackage{a4wide}
\usepackage{pgf,tikz}
\usetikzlibrary{arrows}
\usepackage{epsfig}
\usepackage{color}
\usepackage{mathdots}

\usepackage{hyperref}

\usepackage{authblk}

\usepackage{subcaption}

\newcommand{\beq}{\begin{equation}}  
\newcommand{\eeq}{\end{equation}}  
\newcommand{\bea}{\begin{eqnarray}} 
\newcommand{\eea}{\end{eqnarray}}   
\newcommand{\bear}{\begin{array}}  
\newcommand{\eear}{\end{array}}

\newtheorem{thm}{Theorem}[section] 

\newtheorem{propn}[thm]{Proposition}

\newtheorem{lem}[thm]{Lemma}

\newtheorem{cor}[thm]{Corollary} 

\newenvironment{prf}{\trivlist \item [\hskip 
\labelsep {\bf Proof:}]\ignorespaces}{\qed \endtrivlist}

\theoremstyle{definition}

\newtheorem{exa}[thm]{Example}
\newtheorem{remark}[thm]{Remark}

\newcommand{\D}{{\mathbb D}}

\newcommand{\Z}{{\mathbb Z}}
\newcommand{\C}{{\mathbb C}}
\newcommand{\R}{{\mathbb R}}

\newcommand{\rx}{\mathrm{x}}
\newcommand{\ry}{\mathrm{y}}

\newcommand{\rd}{\mathrm{d}}

\newcommand\la{{\lambda}}
\newcommand\La{{\Lambda}}

\newcommand\al{{\alpha}}
\newcommand\be{{\beta}}
\newcommand\ze{{\zeta}}

\newcommand\om{{\omega}}

\newcommand\si{{\sigma}}

\newcommand\sh{{\cal S}}

\newcommand\ve{{\varepsilon}}

\newcommand\ha{\hat{\alpha}}
\newcommand\hb{\hat{\beta}}
\newcommand\hc{\hat{\gamma}}

\begin{document}


\title{Casting light on shadow Somos sequences}
\author[1]{Andrew N. W. Hone} 
\affil[1]{School of Mathematics, Statistics \&  Actuarial Science, 
University of Kent, 
Canterbury CT2 7FS, U.K.\footnote{ 
e-mail: A.N.W.Hone@kent.ac.uk}
}

\maketitle

\begin{abstract} 
Recently Ovsienko and Tabachnikov considered extensions of Somos and Gale-Robinson sequences, defined over 
the algebra of dual numbers. Ovsienko used the same idea to construct so-called shadow sequences 
derived from other nonlinear recurrence relations exhibiting the Laurent phenomenon, with the 
original motivation  being the hope that these 
examples should lead to an appropriate notion of a cluster superalgebra, incorporating Grassmann variables. 
Here we present various explicit expressions for the shadow of Somos-4 sequences, and describe the solution of a general
Somos-4 recurrence defined over the $\mathbb{C}$-algebra of dual numbers from several different viewpoints: 
analytic formulae in terms of elliptic functions, linear difference equations,  and  Hankel determinants.  
\end{abstract}

\section{Introduction}

\setcounter{equation}{0}

The standard dual numbers, 
which were introduced by Clifford,  take the form $x+y\ve $, where $x,y$ are a pair of real numbers and $\ve^2=0$. 
Triples of such numbers on the dual unit sphere were employed by Study to describe the space of oriented lines in $\R^3$. 
This geometrical interpretation leads to 
contemporary applications of dual numbers in computer vision, while in an algebraic setting they can be used for automatic differentiation. Dual 
numbers also provide the simplest example of an algebra incorporating Grassmann variables, which encode fermionic fields 
in quantum theory.

Cluster algebras, introduced in \cite{fz1}, are a class of commutative algebras with a distinguished set of generators (cluster 
variables) that are defined recursively by a process called mutation, and arise  
in a  variety of different contexts, including Lie theory, Poisson geometry, Teichm\"uller theory and 
integrable systems (see e.g.  \cite{gsv, hkq, in} and references therein). 
A notable feature of cluster algebras is that they exhibit the Laurent phenomenon: each cluster 
variable is a Laurent polynomial in the variables from an initial set of generators (called a seed), with integer coefficients. 
In tandem with the development of cluster algebras, based  on the Caterpillar Lemma \cite{fz}, 
Fomin and Zelevinsky found a systematic way to prove that the Laurent 
phenomenon holds for a wide variety of nonlinear recurrence relations, including many examples 
originally described by Gale \cite{gale}. One such  example is the Somos-4 recurrence, 
namely 
\beq\label{s4} 
x_{n+4}x_n = \al \, x_{n+3}x_{n+1}+\be\, x_{n+2}^2, 
\eeq 
where $\al,\be$ are coefficients. The Laurent property for (\ref{s4}) means that 
the iterates are Laurent polynomials in a set of four initial values $x_0,x_1,x_2,x_3$ 
with coefficients in $\Z[\al,\be]$, 
that is   
\beq\label{lprop}
x_n\in \Z[\al,\be,x_0^{\pm 1},x_1^{\pm 1},x_2^{\pm 1},x_3^{\pm 1}], \quad \forall n\in\Z.  
\eeq 
In particular this implies that the ``classical'' Somos-4 sequence, defined by taking coefficients
$\al=\be=1$ and fixing all initial data to be 1, consists entirely of integers, beginning 
with 
\beq\label{s4seq}
1,1,1,1,2,3,7,23,59,314,1529,8209,83313,\ldots  
\eeq 
(see \cite{oeis}). It was subsequently shown by Fordy and Marsh that the Somos-4 recurrence is generated 
from a cluster algebra defined by a quiver that has periodicity under a specific sequence of mutations: 
the variables $(x_0,x_1,x_2,x_3)$ can be taken as an initial seed, extended by a pair of frozen variables 
$\al,\be$ that do not mutate \cite{FM}.

Given the relevance of superalgebras and supermanifolds in both geometry and theoretical physics, and the myriad 
ways that cluster algebras interact with these different areas, it is currently of interest to find an appropriate 
notion of a cluster superalgebra, incorporating (anticommuting) Grassmann variables in addition to the usual 
(commuting) cluster variables. 
A first step in this direction was taken by Ovsienko \cite{ov1},  who introduced a type of extended quiver, given 
by a hypergraph obtained by adding extra  odd vertices (associated with Grassmann variables) to the usual 
(even) vertices associated with the cluster variables in a cluster algebra defined by a quiver. This allowed the development 
of various examples, including superfriezes, proposed as superalgebra analogues of Coxeter's frieze patterns \cite{mgot, os}, 
and also versions of Somos-4 and higher order Somos-$k$ (or Gale-Robinson) recurrences defined over the dual 
numbers \cite{ot}.   In particular, in \cite{ov1} and \cite{ot} the following dual number generalization of the 
original Somos-4 recurrence with $\al=\be=1$ was considered: 
\beq\label{s4be}
X_{n+4}X_n =X_{n+3}X_{n+1}+(1+\ve)\, X_{n+2}^2; 
\eeq 
while in \cite{ot}, an alternative dual number version was also mentioned, namely 
\beq\label{s4al}
X_{n+4}X_n =(1+\ve)\, X_{n+3}X_{n+1}+ X_{n+2}^2. 
\eeq 
Both of the latter examples define a sequence of dual numbers $X_n=x_n+ y_n\ve $, given a suitable set of four 
initial values $X_0,X_1,X_2,X_3$. It is proved in \cite{ov1} that, within the setting of mutations of extended quivers 
considered there, 
the Laurent property holds in terms of both even and odd variables. As a consequence, for the recurrences (\ref{s4be}) and 
(\ref{s4al}), if the initial values are given by $x_0=x_1=x_2=x_3=1$ and any four integers $y_0,y_1,y_2,y_3$, 
then the whole sequence $(y_n)$ consists of integers.  

More recently, Ovsienko has considered several other examples of 
nonlinear recurrence relations or birational transformations where each variable $x$ is replaced by a dual number
$X=x+ y\ve $, referring to the  corresponding sequence of $y$ values as the shadow sequence \cite{ov2}. For instance, 
the Cassini relation for the Fibonacci sequence produces the convolution of the sequence with itself as a shadow, 
while certain shadow sequences of the Markov numbers appear to be new. 

In this paper we will take the complex numbers as the ambient field, and work with the 
commutative $\C$-algebra of dual numbers $\D=\D(\C)$ given by 
$$ 
\D=\{ \,x+ y \ve \, | \, x,y\in\C, \, \ve^2=0\,\}, 
$$ 
which is isomorphic to the quotient $\C[t]\,/\left<t^2\right>$. Note that the units in $\D$ are the set of elements 
$$ 
\D^*=\{ \,x+ y \ve \in\D\, | \, x\neq 0\,\}, 
$$ 
with 
\beq\label{recip}
(x+y\ve)^{-1}=x^{-1}(1-x^{-1}y\ve).
\eeq 
Then we consider the general Somos-4 recurrence 
for dual numbers $X_n=x_n+ y_n\ve $, given by 
\beq\label{ds4} 
 X_{n+4}X_n = (\al^{(0)}+ \al^{(1)}\ve ) \, X_{n+3}X_{n+1}+(\be^{(0)}+\be^{(1)}\ve )\, X_{n+2}^2, 
\eeq 
 which includes (\ref{s4be}) and 
(\ref{s4al}) as special cases. 
Iteration of (\ref{ds4}) requires that $X_n\in\D^*$ at each step. If $x_n=0$ at some stage, so that $X_n$ is 
not a unit, then a priori it appears impossible to iterate further, 
but one can still consider sequences in $\D$ that satisfy (\ref{ds4}). 
In fact, zero terms in Somos-4 sequences are generically isolated (see the discussion in \cite{honeswart}). 
Moreover, if 
one starts from initial data $X_0,X_1,X_2,X_3\in\D^*$ then the orbit is defined for all $n\in\Z$, 
because the iterates can be obtained by evaluating Laurent polynomials at these four units, 
due to the Laurent phenomenon (see (\ref{lp})  below).

If we separate the above equation into even/odd components (i.e.\ the two components 
of the dual number on each side of the equation), then we obtain the triangular system 
\beq\label{s4even} 
x_{n+4}x_n = \al^{(0)} \, x_{n+3}x_{n+1}+\be^{(0)}\, x_{n+2}^2, 
\eeq 
\small
\beq\label{s4odd}
x_n y_{n+4} - \al^{(0)} x_{n+1}y_{n+3}-2\be^{(0)}x_{n+2}y_{n+2}
 - \al^{(0)} x_{n+3}y_{n+1} +x_{n+4}y_n =  \al^{(1)} x_{n+1}x_{n+3}  +\be^{(1)} x_{n+2}^2. 
\eeq 
\normalsize
The even equation (\ref{s4even})  is just the ordinary Somos-4 recurrence for $x_n$,  with 
coefficients $ \al^{(0)},\be^{(0)}$, while the odd equation  (\ref{s4odd}) is an 
inhomogeneous linear equation for $y_n$, where the coefficients and the inhomogeneity (that is, 
the source term on the right-hand side) are given in terms of $x_n$. In the special case that 
$  \al^{(1)} =\be^{(1)}=0$ and the right-hand side of (\ref{s4odd}) vanishes, the resulting 
homogeneous equation is just the linearization of (\ref{s4even}), corresponding to shadow sequences 
in the sense of \cite{ov2}. 

The original methods for proving the Laurent property for (\ref{s4}), such as the one described by Gale in \cite{gale}, 
or the approach taken by Fomin and Zelevinsky in \cite{fz}, treat the initial data as formal variables but do not 
rely on the auxiliary structure of a cluster algebra. Since $\D$ is a commutative ring, these methods carry over directly to 
(\ref{ds4}), 
immediately yielding the assertion that 
\beq\label{lp}
X_n\in\Z[\al,\be,X_0^{\pm 1}, X_1^{\pm 1}, X_2^{\pm 1}, X_3^{\pm 1}]
\eeq  
for all $n$, where $\al=\al^{(0)}+ \al^{(1)}\ve$, $\be=\be^{(0)}+ \be^{(1)}\ve$. Thus taking even/odd components and using (\ref{recip}), this allows us to state the Laurent 
property for the system (\ref{s4even}),  (\ref{s4odd})  in the following form. 

\begin{lem} 
The system  (\ref{s4even}), (\ref{s4odd}) has the Laurent property, in the sense that 
$$
x_n\in \Z[\al^{(0)},\be^{(0)},x_0^{\pm 1},x_1^{\pm 1},x_2^{\pm 1},x_3^{\pm 1}], 
$$ 
and 
$$ 
y_n\in \Z[\al^{(0)},\al^{(1)},\be^{(0)},\be^{(1)},x_0^{\pm 1},x_1^{\pm 1},x_2^{\pm 1},x_3^{\pm 1},
y_0,y_1,y_2,y_3],  
$$ 
for all $n\in\Z$. 
\end{lem}

The rest of the paper is devoted to providing various representations for 
explicit solutions of (\ref{ds4}), in three different formats: analytic formulae, based on the results of \cite{hones4}; 
an elementary algebraic description, based on the theory of linear difference equations; and Hankel determinant 
expressions, using more recent results on Jacobi continued fractions in \cite{honecf}. 
Each of the next three sections will deal with one 
of these representations, before we end with some conclusions.

\section{Analytic formulae} 

 \setcounter{equation}{0}

To begin with, we paraphrase the main result of \cite{hones4}, and briefly describe it in a way that incorporates certain 
later observations made in \cite{hones5} and \cite{honeswart}. 

\begin{thm}\label{sigma} 
The general solution of the initial value problem for (\ref{s4}) over $\C$ is 
\beq\label{sig}
x_n = A\, B^n \, \frac{\si (z_0+nz)}{\si(z)^{n^2}},  
\eeq 
where $\si(z)=\si(z;g_2,g_3)$ is the Weierstrass sigma function 
for an associated elliptic curve 
\beq\label{weier} 
\ry^2=4\rx^3-g_2\rx-g_3, 
\eeq 
where $A,B\in\C^*$ and $z,z_0,g_2,g_3\in\C$ are explicitly determined by four non-zero initial values and coefficients
$x_0,x_1,x_2,x_3,\al,\be$. 
\end{thm}   

More precisely, 
the solution of the initial value problem for (\ref{s4})  is achieved by considering the sequence of 
ratios
\beq\label{dn} 
d_n = \frac{x_{n+1}x_{n-1}}{x_n^2}, 
\eeq 
which satisfy a recurrence of  second order, namely  
\beq\label{qrt} 
d_{n+1}d_{n-1} = \frac{\al d_n+\be}{d_n^2}. 
\eeq 
The latter recurrence  is equivalent to a birational map of the plane, being an example of a symmetric QRT map 
\cite{qrt}, and it has a rational conserved quantity 
\beq \label{J}
J=d_{n}d_{n-1}+ \al \left(\frac{1}{d_n}+\frac{1}{d_{n-1}}\right) + \frac{\be}{d_nd_{n-1}}
\eeq 
that defines a pencil of biquadratic plane curves, each of which is generically of genus 1 (except for a finite 
number of singular fibres, for certain special values of $J$ \cite{duistermaat}) and hence isomorphic to a 
 cubic curve of the form (\ref{weier}). 
Explicit computations with 
Weierstrass functions then show that the formula (\ref{sig}) does indeed provide a solution of 
the Somos-4 recurrence (\ref{s4}), with the coefficients being given by 
\beq\label{ab} 
\al =\frac{\si (2z)^2}{\si (z)^8}, \qquad \be = -\frac{\si(3z)}{\si(z)^9}, 
\eeq 
allowing $g_2,g_3$ to be obtained in terms of $\al,\be,J$ via the expressions 
\beq\label{g23}
g_2 =12\la^2-2J, \quad g_3=4\la^3-g_2\la-\al, \quad \la=\frac{1}{3\al}\left(\frac{J^2}{4}-\be\right), 
\eeq 
and then (up to fixing an overall sign) $z,z_0$ are determined modulo periods by the elliptic integrals
\beq\label{eint}
z=\int^\la_\infty \frac{\rd\rx}{\ry}, \qquad z_0=\int^{\la-d_0}_\infty \frac{\rd\rx}{\ry}, 
\eeq  
so that once these are given the values of  $A,B$ can be found from the initial data.  
(Note that when $g_2^3-27g_3^2=0$, corresponding to a singular fibre, the formula (\ref{sig}) is still
valid in terms of suitable trigonometric/rational limits of the sigma function.)

For the purposes of what follows, it is important to note that from (\ref{dn}) the 
expression (\ref{J}) is equivalent to a 
rational conserved quantity for the Somos-4 recurrence (\ref{s4}) itself, that is 
\beq\label{Jx} 
J = \frac{x_n^2x_{n+3}^2 +\al\, (x_{n+1}^3x_{n+3}+x_nx_{n+2}^3) +\be\,x_{n+1}^2x_{n+2}^2 }
{x_nx_{n+1}x_{n+2}x_{n+3}}. 
\eeq

An elementary observation, which is the basis for the application of dual numbers to automatic differentiation, 
is that for any differentiable function $\Phi$ and $X=x+y\ve\in\D$, the value of $\Phi$ at $X$ is 
$$ 
\Phi(X)=\Phi(x)+\Phi'(x)y\ve.
$$ 
We can use this result to describe solutions of (\ref{ds4}) in analytic form.

\begin{propn}\label{anal}
Given parameters $A=A^{(0)}+A^{(1)}\ve$, 
$B=B^{(0)}+B^{(1)}\ve$, 
$G_2=g_2+\tilde{g}_2\ve$, 
$G_3=g_3+\tilde{g}_3\ve$, $Z_0=z_0+\tilde{z}_0\ve$, 
$Z=z+\tilde{z}\ve\in\D$, 
the expression 
\beq\label{Xsig}
X_n = A\, B^n \, \frac{\si (Z_0+nZ)}{\si(Z)^{n^2}}  
\eeq 
satisfies 
the dual Somos-4 recurrence (\ref{ds4}) with coefficients 
$$ 
\al =\al^{(0)}+ \al^{(1)}\ve=\frac{\si (2Z)^2}{\si (Z)^8}, \quad 
\be =\be^{(0)}+ \be^{(1)}\ve = -\frac{\si(3Z)}{\si(Z)^9}, 
$$
provided that $A,B\in\D^*$ and $\si(Z)=\si(Z;G_2,G_3)\in\D^*$. 
In terms of even/odd components, this can be written as 
\beq\label{xyan}
X_n =x_n +x_n\left(\frac{A^{(1)}}{A^{(0)}}+ \frac{B^{(1)}}{B^{(0)}} \, n
+\tilde{z}_0\zeta(z_0+nz) + (\tilde{z}\partial_z+\tilde{g}_2\partial_{g_2}+\tilde{g}_3\partial_{g_3})\log x_n
\right)\,\ve, 
\eeq
where $x_n$ denotes the right-hand side of (\ref{sig}) with the replacement 
$A\to A^{(0)}$, $B\to B^{(0)}$, and the other parameters $z,z_0,g_2,g_3$ being the same, 
$\zeta(z)=\zeta(z;g_2,g_3)$ is the Weierstrass zeta function, and $\partial$ denotes a partial 
derivative.
\end{propn}

\begin{prf} As pointed out in \cite{hones5}, the fact that the analytic expression (\ref{sig}) satisfies 
a Somos-4 relation with coefficients given by (\ref{ab}) is a direct consequence of the three-term 
identity for the sigma function, which can be viewed as an algebraic identity between power series. 
The sigma function $\si(z;g_2,g_3)$ is a holomorphic function of the argument $z$, and is also 
holomorphic in the two invariants $g_2,g_3$. Therefore the same proof carries over to the commutative algebra 
$\D$, with the requirement that $\si(Z)\in\D^*$ so that one can divide by suitable powers of this quantity, 
and $A,B\in\D^*$ so that the even part $x_n\neq 0$.  The solution 
depends analytically on six dual parameters, so it can be expanded as $X_n=x_n+y_n\ve$ by differentiating 
with respect to each of these parameters, giving 
$$ 
y_n =(A^{(1)}\partial_{A^{(0)}}+ B^{(1)}\partial_{B^{(0)}}+ 
+\tilde{z}_0\partial_{z_0}+\tilde{z}\partial_{z}
+\tilde{g}_2\partial_{g_2}+\tilde{g}_3\partial_{g_3})\,x_n.
$$
It is convenient to rewrite this in terms of logarithmic derivatives of $x_n$, so that 
$y_n$ is written as $x_n$ times a sum of such derivatives, and the first three derivatives are 
straightforward to calculate explicitly, as in (\ref{xyan}).
\end{prf}

We believe that the analytic formula (\ref{Xsig}) represents the general solution of the dual Somos-4 recurrence (\ref{ds4}).
However, a complete proof would require showing that, given $\al,\be\in\D$, one can solve the initial value problem for (\ref{ds4}) by determining the four dual parameters $G_2,G_3,Z_0,Z$. The difficulty lies in constructing a suitable 
theory of elliptic functions and uniformization of elliptic curves over the dual numbers, which in 
particular would provide analogues of the elliptic integrals (\ref{eint}). An abstract theory of super Riemann surfaces 
is available \cite{witten}, but to the best of our knowledge the genus 1 case has not been developed 
explicitly in terms of Weierstrass functions with odd components. Rather than trying to pursue this further here, 
we will take the formula (\ref{xyan}) as the starting point for a more straightforward algebraic approach in the next 
section.

\section{Solution of linear difference equation} 

\setcounter{equation}{0}

Since the system consisting of (\ref{s4even}) and (\ref{s4odd}) is triangular, and we already have the 
general solution of (\ref{s4even}), given by Theorem \ref{sigma}, the main outstanding problem is to understand the 
solutions of (\ref{s4odd}). We beginning by considering the special  case $\al^{(1)} =\be^{(1)}=0$, 
which is just a homogeneous linear 
equation for $y_n$, namely 
\beq\label{shad}
x_n y_{n+4} - \al^{(0)} x_{n+1}y_{n+3}-2\be^{(0)}x_{n+2}y_{n+2}
 - \al^{(0)} x_{n+3}y_{n+1} +x_{n+4}y_n = 0.  
\eeq 
The solutions of the latter equation are the shadow Somos-4 sequences, in the sense of \cite{ov2}. 

However, the equation (\ref{shad}) is just the linearization of the Somos-4 recurrence 
(\ref{s4even}), and we have the explicit solution of this, which (for fixed coefficients $\al^{(0)}$, $\be^{(0)}$) 
depends on four arbitrary parameters; these can be taken to be  
$A^{(0)}$, $B^{(0)}$, $z_0$ and $J^{(0)}$, where the last is the value of the first integral 
(\ref{Jx}). Thus, by adapting a technique that goes back to Lie's work on ordinary differential equations, 
we can observe that the derivative of the solution $x_n$ with respect to each of these parameters  provides 
a solution of the linearized equation, and thus we obtain four linearly independent solutions, which span 
the whole vector space of solutions of (\ref{shad}).

\begin{lem}\label{lin}
The Somos-4 shadow equation  (\ref{s4odd}) has four linearly independent solutions, given by 
$y_n^{(i)}=x_n$, $y_n^{(ii)}=nx_n$, $y_n^{(iii)}=x_n\zeta(z_0+nz)$ and 
\beq\label{y4}
y_n^{(iv)}=x_n\partial_{J^{(0)}}\log x_n 
=x_n\left(\frac{\rd z}{\rd J^{(0)}}\partial_z 
+\frac{\rd g_2}{\rd J^{(0)}}\partial_{g_2} 
+ +\frac{\rd g_3}{\rd J^{(0)}}\partial_{g_3} \right)\log x_n .
\eeq 
\end{lem} 

\begin{prf} 
The derivative of the formula (\ref{sig}) with respect to the scaling parameter $A$ is just proportional to $x_n$, 
while the derivative with respect to $B$ is proportional to $nx_n$. Differentiating the numerator in (\ref{sig}) 
with respect to 
$z_0$ produces $x_n$ times the Weierstrass zeta function, but the dependence on $J^{(0)}$ is more complicated,  
as both $z$ and the invariants $g_2,g_3$ vary with this parameter.
\end{prf} 

\begin{remark} The fact that $y_n=x_n$ is a solution of (\ref{s4odd}) was already noted in \cite{ot}.
\end{remark}

The parameters $A,B$ in (\ref{sig}) correspond to the action of a two-parameter scaling group, 
$x_n\to AB^n\,x_n$ for $(A,B)\in(\C^*)^2$, which leaves the Somos-4 recurrence invariant. The freedom to vary 
$z_0$ corresponds to an arbitrary choice of initial point on the elliptic curve (\ref{weier}), viewed as a 
complex torus $\C/\Lambda$ where $\Lambda$ is the period lattice. 
The first three terms in the odd part of (\ref{xyan}) are clearly proportional to $y_n^{(i)}$, $y_n^{(ii)}$ and 
$y_n^{(iii)}$.
We will shortly present an algebraic method to calculate the third independent shadow solution $y_n^{(iii)}$. 

As 
for the fourth solution  $y_n^{(iv)}$, it is much more complicated, because it involves the modular derivatives  
$\partial_{g_2} $ and $\partial_{g_3}$. Indeed, if we fix the coefficients $\al^{(0)}$, $\be^{(0)}$, 
then we can consider the following system of three equations in terms of the Weierstrass 
$\wp$ function and its derivatives: 
\beq\label{abjsys} 
\al^{(0)}=\wp'(z)^2, \qquad \frac{\be^{(0)}}{\al^{(0)}}=\wp(2z)-\wp(z), 
\qquad J^{(0)}=\wp''(z); 
\eeq 
the first two equations are just equivalent to (\ref{ab}), corresponding to identities between elliptic functions 
(see \cite{hones4, hones5}). Differentiating each of these equations with respect to $J^{(0)}$ 
and applying the chain rule gives a system of three linear equations for 
the derivatives of $z,g_2,g_3$ with respect to $J^{(0)}$, appearing in (\ref{y4}). 
As for  the modular derivatives of the sigma function, note that it 
is a weighted homogeneous function 
of $z,g_2,g_3$, so it satisfies the linear partial differential equation (PDE)
$$ 
(4g_2\partial_{g_2}+6g_3\partial_{g_3}  -z\partial_z+1)\sigma(z;g_2,g_3)=0. 
$$
As was shown by Weierstrass \cite{weierstrass}, it satisfies another linear PDE which is first order 
in $g_2,g_3$ and second order in $z$: this is equivalent to the well-known heat equation for the Jacobi theta function. 
(It also satisfies a fourth order nonlinear ODE in $z$ that can be written in Hirota bilinear form \cite{ee}, 
which can be viewed as the continuous analogue of Somos-4.)  
In principle one could combine all these facts to rewrite (\ref{y4}) in terms of $z$ derivatives only. However, 
we do  not need to calculate these analytic expressions in more detail, because in due course we will 
present an algebraic way to characterize a fourth independent shadow solution, which is much more 
straightforward. 

In order to derive the third shadow solution in a more algebraic way, we consider the following coupled pair of recurrence 
relations with parameters $u,f$, which was considered in \cite{honecf} as 
the simplest example of a map arising from Jacobi continued fractions in hyperelliptic function fields, 
based on a construction introduced by van der Poorten \cite{vdp}: 
\beq\label{dtoda}
\begin{array}{rcl} 
v_n & = & -v_{n-1}+u/d_n, \\ 
d_{n+1} & = & -d_n-v_n^2-f. 
\end{array} 
\eeq 
The above system defines an integrable map, and here we will present its analytic solution. 
\begin{propn} 
The map $(v_{n-1},d_n)\mapsto (v_n,d_{n+1})$ defined by (\ref{dtoda}) is 
integrable: it preserves the symplectic form $\rd v_{n-1}\wedge \rd d_n$, and has the first integral 
$$H=d_n(v_{n-1}^2+d_n+f)-uv_{n-1}.$$
The general solution of the map can be written in terms of elliptic functions, 
up to an overall choice of sign, as 
\beq\label{vdan}
d_n=\wp (z)-\wp(z_0+nz), \qquad v_n=\pm \big(\ze (z_0+(n+1)z)- \ze (z_0+nz)-\ze(z)\big), 
\eeq 
with the parameters given by 
$$ 
u=\pm \wp'(z), \qquad f=-3\wp(z),
$$ 
on the level set $H=-\frac{1}{2}\wp''(z)$. 
\end{propn}
\begin{prf}
It is straightforward to verify from the map (\ref{dtoda}) that 
$\om = \rd v_{n-1}\wedge \rd d_n=\rd v_{n}\wedge \rd d_{n+1}$, so it is a symplectic map, and 
a short calculation also shows that $H$ is independent of $n$, so this is a discrete 
integrable system with one degree of freedom.  
By Proposition 5.1 of \cite{honecf}, on each level set $H=\,$const, the sequence of 
values of $d_n$ coincides with an orbit of (\ref{qrt}), if the coefficients  and 
the value of the first integral (\ref{J}) are identified as 
\beq\label{abjform} 
\al=u^2, \qquad \be=u^2(v^2+f), \qquad 
J=2uv=-2H.
\eeq
 As for the analytic solution of the map, it was shown in \cite{hones4} 
that  when $x_n$ is given by (\ref{sig}), the corresponding solution of (\ref{qrt}), associated 
to it via (\ref{dn}), is given in terms of the Weierstrass $\wp$ function by the first formula in 
(\ref{vdan}), and when this expression for $d_n$ is substituted into the second component of 
(\ref{dtoda}), the explicit formula for $v_n$ is found  from taking the square root 
in the left-hand side of the elliptic function identity 
$$\big(\ze(a)+\ze(b)+\ze(c)\big)^2=\wp(a)+\wp(b)+\wp(c), 
\quad\mathrm{for}\quad a+b+c\equiv 0 \bmod\La.
$$
Then given the two formulae in (\ref{vdan}), another identity between elliptic functions 
verifies that the first component of the system (\ref{dtoda}) is satisfied for all $n\in\Z$, 
with the parameter $u=\pm \sqrt{\al}=\pm \wp'(z)$, making the same choice of sign as for 
$v_n$. 
\end{prf}

\begin{remark} Given the canonical Poisson bracket $\{\,d_n,v_{n-1}\,\}=1$ 
associated with $\om$, the first 
integral $H$ defines the Hamiltonian vector field $\{ \,\cdot, H\,\}$, whose flow  commutes with the 
shift $n\to n+1$ corresponding to the map. A direct calculation shows that, for a fixed choice of scale, this flow implies 
that the sequence of $d_n$ and $v_n$ satisfy the set of differential equations 
\beq\label{toda}
\begin{array}{rcl} 
\dot{v}_n & = & d_n-d_{n+1}, \\ 
\dot{d}_n  & = & d_n(v_{n-1}-v_n),  
\end{array} 
\eeq 
for all $n\in\Z$, which (up to squaring one of the variables) is just the infinite Toda lattice written in 
Flaschka coordinates. Thus the simultaneous solutions of the map and the Hamiltonian flow provide 
genus 1 solutions of the Toda lattice, and for genus $g>1$ an analogous statement holds for the other maps 
constructed in  \cite{honecf}; further details will be presented elsewhere. 
\end{remark}

\begin{cor}\label{y3} 
Up to subracting off multiples of $y_n^{(i)}=x_n$ and $y_n^{(ii)}=nx_n$ and overall scale, 
the third independent shadow solution of Somos-4 can be obtained from the 
associated solution of the map (\ref{dtoda}), in the form  
\beq\label{vsum}y_n^{(iii)}=-x_n\sum_{j=0}^{n-1}v_j, \qquad n\geq 0\eeq 
(where the case $n=0$ corresponds to an empty sum).  
\end{cor}
\begin{prf}
If we take the plus sign in the formula for $v_n$ in (\ref{vdan}), then 
we have a telescopic sum 
$$-x_n\sum_{j=0}^{n-1}v_j =x_n\big(n\ze(z)+\ze(z_0)-\ze(z_0+nz)\big)$$ for $n\geq 0$, and after subtracting off 
$\ze(z_0) y_n^{(i)}$ and $\ze(z)y_n^{(ii)}$ this is just the third independent solution $y_n^{(iii)}$ in Lemma \ref{lin}, 
up to an overall sign. 
\end{prf}

We now present an observation that dramatically simplifies the solution of the dual Somos-4 equation, and allows 
the order of (\ref{s4odd}) to be reduced from four to three. The point is that, since $\D$ is a commutative algebra, the calculation which shows that (\ref{Jx}) is a conserved quantity for (\ref{s4}) carries over directly to the dual numbers, so the same 
expression with the replacement $x_n\to X_n$ provides a first integral  $J=J^{(0)}+J^{(1)}\ve\in\D$ for  (\ref{ds4}). 
This implies that the system of even/odd equations has a pair of rational first integrals, namely $J^{(0)}$, which 
is just the first integral for the even part (\ref{s4even}) given by the original expression (\ref{Jx}) in terms of 
$x_n$ alone, but with parameters $\al\to\al^{(0)}$, $\be\to\be^{(0)}$,  and $J^{(1)}$, which depends on both sets of variables $x_n$ and $y_n$, and is linear in $y_n$. 
The most efficient way to calculate $J^{(1)}$ is to clear the denominator in the formula for $J\in\D$, multiplying both 
sides of the relation by $X_nX_{n+1}X_{n+2}X_{n+3}$, to obtain a polynomial relation, and then the leading order (even) part of the 
resulting expression just returns the usual formula for $J^{(0)}$, while the $O(\ve)$ (odd) part gives a linear equation in 
 $J^{(1)}$, but also contains $J^{(0)}$. Upon eliminating $J^{(0)}$ from the latter expression, the desired formula 
for $J^{(1)}$ is obtained.

\begin{lem} \label{jlem} 
In addition to the conserved quantity $J^{(0)}$, given by replacing the 
parameters  $\al\to\al^{(0)}$, $\be\to\be^{(0)}$ in the right-hand side of (\ref{Jx}),  the system 
consisting of  (\ref{s4even}) and  (\ref{s4odd}) has a second rational first integral, namely 
 \beq\label{j1form} 
J^{(1)} = \frac{D_n-\sum_{j=0}^3C^{(j)}_n x_{n+j}^{-1} y_{n+j}}
{x_nx_{n+1}x_{n+2}x_{n+3}},
\eeq
where 
$$ \begin{array}{rcl}
C_n^{(0)} & = & \al^{(0)}x_{n+1}^3x_{n+3}+\be^{(0)}x_{n+1}^2x_{n+2}^2-x_n^2x_{n+3}^2, \\
C_n^{(1)} & = & \al^{(0)}x_nx_{n+2}^3-2\al^{(0)}x_{n+1}^3x_{n+3}-\be^{(0)}x_{n+1}^2x_{n+2}^2+x_n^2x_{n+3}^2, \\
C_n^{(2)} & = &  -2\al^{(0)}x_nx_{n+2}^3+\al^{(0)}x_{n+1}^3x_{n+3}-\be^{(0)}x_{n+1}^2x_{n+2}^2+x_n^2x_{n+3}^2, \\
C_n^{(3)} & = &  \al^{(0)}x_nx_{n+2}^3+\be^{(0)}x_{n+1}^2x_{n+2}^2-x_n^2x_{n+3}^2, \\
D_n & = &  \al^{(1)}x_nx_{n+2}^3+\al^{(1)}x_{n+1}^3x_{n+3}+\be^{(1)}x_{n+1}^2x_{n+2}^2. \\
\end{array}
$$
\end{lem}

For fixed $J^{(1)}$, the equation (\ref{j1form}) can be rearranged as a linear inhomogeneous difference equation 
of third order for $y_n$, that is 
\beq\label{lin3}
L_n (y_n) = F_n, 
\eeq 
where $L_n$ is the linear difference operator 
\beq\label{lop}
L_n = C_n^{(3)} x_{n+3}^{-1} \sh^3
+ C_n^{(2)} x_{n+2}^{-1} \sh^2
+ C_n^{(1)} x_{n+1}^{-1} \sh
+ C_n^{(0)} x_{n}^{-1}, 
\eeq  
given in terms of the shift operator $\sh$ that sends $n\to n+1$, and 
\beq\label{inhg}
F_n = D_n -J^{(1)}x_nx_{n+1}x_{n+2}x_{n+3}.
\eeq
Thus we have succeeded in reducing the order of (\ref{s4odd}) by one, as claimed. 
The corresponding third order homogeneous equation which arises when 
$\al^{(1)}=\be^{(1)}=J^{(1)}=0$, namely 
\beq\label{hgs} 
L_n (y_n) =0,
\eeq 
is nothing other than the linearization of the equation defining $J^{(0)}$, and  three linearly independent solutions 
are obtained by varying the analytic solution (\ref{sig}) with respect to the three parameters that do not depend on this modular 
quantity: in other words, up to rescaling  and taking linear combinations, they are given by the 
first three shadow solutions $y^{(i)}_n,y^{(ii)}_n,y^{(iii)}_n$ in Lemma \ref{lin}.    
Then a fourth independent Somos-4 shadow solution is found by taking $\al^{(1)}=\be^{(1)}=0$ but $J^{(1)}\neq0$. 

The form of the solutions of the  homogeneous equation (\ref{hgs}) is made more 
transparent by setting $y_n=x_nY_n$. This gives 
$$L_n(y_n)=x_nx_{n+1}x_{n+2}x_{n+3}\tilde{L}_n (Y_n)=0,$$ where 
$$ \tilde{L}_n =  \left(\Big(\frac{\be^{(0)}}{d_{n+1}d_{n+2}}-d_{n+1}d_{n+2}\Big)(\sh+1)
 +\frac{\al^{(0)}}{d_{n+2}} \, \sh  
+\frac{\al^{(0)}}{d_{n+1}}
\right)(\sh -1)^2.
$$
So clearly $Y_n=1$ and $Y_n=n$ lie in the kernel of the latter operator, corresponding to 
 $y^{(i)}_n$ and $y^{(ii)}_n$.

Finally, we can construct the general solution of (\ref{s4odd}) by applying the discrete analogue of the method of variation 
of parameters to find an arbitrary element in the affine space of solutions of (\ref{lin3}). In other words, we can write the 
solution as the sum of a particular integral plus an arbitrary linear combination of solutions of the  
homogeneous equation (\ref{hgs}). 
For variation of parameters, the initial ansatz is to write the solution of (\ref{lin3}) in the form 
\beq\label{vopsol}
y_n = \sum_j f_n^{(j)} y_n^{(j)}, 
\eeq
where the index $j$ ranges over the three lower case Roman numerals $i,ii,iii$, and then impose the 
constraints that 
\beq\label{req} 
 \sum_j (f_{n+1}^{(j)}-f_n^{(j)}) y_{n+1}^{(j)}=0= \sum_j (f_{n+1}^{(j)}-f_n^{(j)}) y_{n+2}^{(j)},
\eeq 
which together imply that $y_{n+1}= \sum_j f_n^{(j)} y_{n+1}^{(j)}$, 
 $y_{n+2}= \sum_j f_n^{(j)} y_{n+2}^{(j)}$. 
Putting all this into   (\ref{lin3}) gives 
$$ 
L_n(y_n) = C_n^{(3)}x_{n+3}^{-1} \sum_j (f_{n+1}^{(j)}-f_n^{(j)}) y_{n+3}^{(j)} 
+ \sum_j f_n^{(j)}L_n ( y_{n}^{(j)})=C_n^{(3)}x_{n+3}^{-1} \sum_j (f_{n+1}^{(j)}-f_n^{(j)}) y_{n+3}^{(j)} , 
$$ 
which must equal $F_n$. 
Combining the latter equality with the two constraints (\ref{req}) gives a linear system for the 
three differences   $f_{n+1}^{(j)}-f_n^{(j)}$, which is solved to yield 
\begin{thm}\label{gensol} 
The general solution of (\ref{lin3}) can be written in the form  (\ref{vopsol}), where 
\beq\label{varp}
f_n^{(j)}=f_0^{(j)}+\sum_{k=0}^{n-1}v_k^{(j)}, \qquad j=i,ii,iii,  \qquad \mathrm{for}\,\,n\geq 0, 
\eeq
with
\beq\label{vop}
\left(\begin{array}{c} v_n^{(i)} \\ v_n^{(ii)} \\ v_n^{(iii)} \end{array}
\right) 
=
\frac{x_{n+3}F_n}{C_n^{(3)}}
\left|\begin{array}{ccc} y_{n+1}^{(i)}  & y_{n+1}^{(ii)} &  y_{n+1}^{(iii)} \\ 
y_{n+2}^{(i)}  & y_{n+2}^{(ii)} &  y_{n+2}^{(iii)} \\ 
y_{n+3}^{(i)}  & y_{n+3}^{(ii)} &  y_{n+3}^{(iii)} 
\end{array}
\right|^{-1}   
\left(\begin{array}{c} y_{n+1}^{(ii)}  y_{n+2}^{(iii)}-y_{n+1}^{(iii)}  y_{n+2}^{(ii)}  \\  
y_{n+1}^{(iii)}  y_{n+2}^{(i)}-y_{n+1}^{(i)}  y_{n+2}^{(iii)}\\  
y_{n+1}^{(i)}  y_{n+2}^{(ii)}-y_{n+1}^{(ii)}  y_{n+2}^{(i)}
\end{array} \right) 
.
\eeq 
The 
three  parameters $f_0^{(j)}$ for $j=i,ii,iii$ in (\ref{varp}) are arbitrary, 
and together with the freedom to  choose $J^{(1)}$  arbitrarily, this in turn 
provides the general solution of (\ref{s4odd}). 
\end{thm}

\begin{table}[h!]
  \begin{center}
    \caption{Four independent shadows of the original Somos-4 sequence (\ref{s4seq}).}
    \label{factortable} 
\scalebox{0.9}{
    \begin{tabular}{ | r|| r| r| r| r |r|r|r|r|r|r|r|r|r|r|} %
\hline
      $n$ & -1 & 0 &1&2&3&4&5&6&7&8&9&10 & 11 & 12  \\
\hline 
&&&&&&&&&&&&&&\\
$y_n^{(i)}$  & 1 & 1 & 1& 1 & 2 & 3 & 7 & 23 & 59 & 314 & 1529 & 8209 & 83313 & 620297 
 \\ 
&&&&&&&&&&&&&&\\
$y_n^{(ii)}$  & -1 & 0 & 1 & 2 & 6 & 12 & 35 & 138 & 413 & 2512 & 13761 & 82090 & 916443 & 7443564 
\\ 
&&&&&&&&&&&&&&\\
 $y_n^{(iii)}$  & 0 & 0 & 1 & 1 & 3 &  7 & 15 & 70 & 202 & 1107 & 6906 & 36386 & 420371 & 3594979 
 \\ 
&&&&&&&&&&&&&&\\
$y_n^{(iv)}$  & 0 & 0 & 0 & 1 & 1 & 3 & 10 & 22 & 108 & 472 & 2174 & 17792 & 120536 & 1161627
 \\
&&&&&&&&&&&&&&\\
\hline 
    \end{tabular}
}
  \end{center}
\end{table}

\begin{exa}\label{orig} As a particular example of solving (\ref{lin3}), let us consider the original Somos-4 sequence 
(\ref{s4seq}), and obtain four independent shadow sequences. For ease of comparison with Example 4.2 in \cite{honecf}, 
it is convenient to index the original sequence so that $x_{-1}=x_0=x_1=x_2=1$, and 
then the first shadow sequence $y_n^{(i)}$ is just given by the same sequence, starting with index $n=-1$, extending to 
a sequence of 
positive terms for all $n\in\Z$ (as due to a symmetry it repeats the same values when run in reverse), while a second 
independent sequence is $y_n^{(ii)}=nx_n$, starting with $y_{-1}^{(ii)}=-1$,  $y_{0}^{(ii)}=0$,   $y_{1}^{(ii)}=1$,  $y_{2}^{(ii)}=2$, 
which is positive for all $n\geq 1$. As for a third independent sequence, it is found by applying Corollary \ref{y3},  
using the system (\ref{dtoda}) 
with parameters $u=-1$, $f=-3$ 
to generate a sequence of pairs $(v_n,d_{n+1})$ starting from $v_0=-1$, $d_1=1$ (and note that we also have $d_0=1$ and $H=-J/2=-2$ in this case). 
Then from (\ref{vsum}) we find $y_{0}^{(iii)}=0$,   $y_{1}^{(iii)}=1$,  $y_{2}^{(iii)}=1$, and going one step back  with the 
homogeneous equation (\ref{hgs}) shows that
$y_{-1}^{(iii)}=0$. From the first equation in (\ref{dtoda}) with $u=-1$ it follows that $-(v_n+v_{n-1})=1/d_n>0$ for all $n$, and so 
from (\ref{vsum}) and the initial value $v_0=-1$ this implies that this third shadow sequence is positive 
whenever $n\geq 1$. 
Finally, for a fourth independent shadow sequence we must solve   (\ref{lin3}) 
with $\al^{(1)}=\be^{(1)}=0$ and a non-zero value of  $J^{(1)}$, so we fix $J^{(1)}=-1$ and 
take $y_n^{(iv)}=0$ for $n=-1,0,1$, giving $y_2^{(iv)}=1$. Empirical evidence suggests the conjecture that this fourth sequence should be positive 
for all $n\geq 2$, but we do not have a proof. Table \ref{factortable} presents the first few  values in these shadow sequences.   
\end{exa}

\section{Hankel determinant formulae} 

\setcounter{equation}{0}

It was conjectured by Barry and proved by various authors that certain Somos-4 sequences could be expressed as Hankel determinants \cite{ch, chx, xin}. 
In \cite{honecf} we showed that these results can be unified and further generalized 
by applying  van der Poorten's work on Jacobi continued fractions (J-fractions) in hyperelliptic function fields \cite{vdp}. 
In the genus 1 case, one expands a certain function $G$ on a quartic curve as a J-fraction, that is 
\beq\label{gcla} 
G
=\cfrac{s_0}{X+v_1 -\cfrac {d_2}{ X+v_2-\cfrac{d_3}{X+v_3 -\cdots}} },  
\eeq
where $X^{-1}$ is a local parameter around one of the points at infinity, and the recursion 
relation for the continued fraction leads to the map (\ref{dtoda}) for $(v_{n-1},d_n)$. The numerators and denominators 
of the convergents provide associated orthogonal polynomials, and standard results imply that 
the power series expansion  of the generating function  near $X=\infty$, that is 
$G=\sum_{j\geq 1} s_{j-1}X^{-j}$, provides a sequence of 
moments $(s_j)$ such that $d_n$ and $v_n$ can be written in terms of ratios of the corresponding determinants $\Delta_n$ 
or $\Delta_n^*$ of 
Hankel/bordered Hankel type, respectively.    
In particular, the solution of (\ref{qrt}) is given by $d_n=\Delta_{n-2}\Delta_n/\Delta_{n-1}^2$, and in \cite{honecf} we 
further 
showed how this result extends to the solution of an integrable symplectic map associated with the J-fraction 
expansion of a function on a hyperelliptic curve 
of any genus $g>1$. 

Since $\D$ is a commutative algebra,  identities for Hankel determinants carry over directly to suitable 
sequences of moments $s_j\in\D$, and allow the solution of (\ref{ds4}) to be expressed in the same form as 
for Somos-4 sequences over $\C$. Thus, simply by writing a  dual number version of the statement of 
Theorem 5.2 in \cite{honecf} 
(also making use of the formulae in Theorem 4.1 therein for the particular case of genus 1), 
and taking care that certain parameters should be units, we arrive at the following result. 
 
\begin{thm}\label{hank}
Given arbitrary $\ha,\hb,\hc,s_0\in\D^*$ and $s_1\in\D$, define a sequence of dual numbers $(s_j)_{j\geq 0}$ by
the recursion  
\beq\label{srec}
s_j = \ha s_{j-2} +\hb\sum_{i=0}^{j-2}s_is_{j-2-i}+\hc \sum_{i=0}^{j-3}s_is_{j-3-i}, \qquad j\geq 2, 
\eeq 
and form the associated sequence of Hankel determinants 
\beq\label{hankel} 
\Delta_n = 
\left| \begin{array}{cccc} 
s_0     & s_1 & \cdots & s_{n-1}   \\ 
s_1     &         & \iddots &        \vdots       \\
\vdots & \iddots &          &           \vdots    \\        
s_{n-1}     &  \cdots & \cdots & s_{2n-2}  
\end{array} 
\right| =\det (s_{i+j-2})_{i,j=1,\ldots, n}
\eeq 
for $n\geq 1$, with the usual convention that $\Delta_0=1$. Then 
$$ 
X_n = \Delta_{n-1} \qquad \mathrm{for}\,\, n\geq 1
$$ 
is a solution of the dual Somos-4 recurrence (\ref{ds4}) with coefficients given by  
$$ 
\al = U^2, \qquad \be =\al F + \frac{1}{4}J^2, 
$$ 
where 
$$ 
U=-s_0\hc-s_1\hb, \qquad F=-\ha-2s_0\hb, 
$$ 
and the value of the first integral is fixed to be 
$$ 
J=2(s_0\ha\hb +s_0^2\hb^2+s_1 \hc)
.
$$
\end{thm}

\begin{remark} 
The Hankel determinant expression above only depends on 5 dual number parameters, which is one less 
than is required to produce the general solution of (\ref{ds4}): there are four initial values and two coefficients 
in the initial value problem for Somos-4. However, the missing parameter can be recovered 
with the scaling symmetry $X_n\to AX_n$ for $A\in\D^*$. The other scaling symmetry $X_n\to B^n X_n$ for 
$B\in\D^*$ just rescales the value of $s_0$, which is an overall multiplier in the moment sequence. 
As discussed in \cite{honecf}, there is another moment sequence which provides Hankel determinant formulae 
for negative indices $n$, and in general it is necessary to apply these two scaling symmetries in order to glue 
the two sequences of Hankel determinants together into a valid Somos-4 sequence for all $n\in\Z$.  
\end{remark}

\begin{exa} To obtain Hankel determinant formulae for  the fourth shadow sequence, as in Example \ref{orig}, we take 
dual numbers that are $O(\ve)$ perturbations of the corresponding quantities in Example 4.2 of \cite{honecf}. 
Setting 
\beq\label{spars}
\ha= 1-4\ve, \quad \hb=1, \quad \hc=1-\tfrac{3}{2}\ve, \quad s_0=1+\ve, \quad s_1=\tfrac{1}{2} \ve
\eeq 
gives $U=-1$, $F=-3+2\ve$, and hence 
$$ 
\al=1, \qquad \be=1, \qquad J=4-\ve, 
$$ 
and this corresponds to $\al^{(0)}=\be^{(0)}=1$, $J^{(0)}=4$ and 
$\al^{(1)}=\be^{(1)}=0$, $J^{(1)}=-1$ as required. 
Then upon iterating the recursion (\ref{srec}) with parameters and initial values as in (\ref{spars}), 
the sequence $(s_j)$ is found to be 
$$ 
1+\ve, \tfrac{1}{2} \ve, 2-\ve, 1+2\ve, 6-6\ve,7+2\ve,24-28\ve,41-23\ve, 115-154\ve, 
236-\tfrac{527}{2}\ve, \ldots , 
$$  
and 
the corresponding sequence of Hankel determinants 
begins with $\Delta_0=1$, $\Delta_1=1+\ve$, 
$$
\Delta_2
=\left|\begin{array}{cc} 1+\ve & \tfrac{1}{2}\ve \\ \tfrac{1}{2}\ve & 2-\ve \end{array}\right|=2+\ve, 
\quad
\Delta_3= \left|\begin{array}{ccc} 1+\ve & \tfrac{1}{2}\ve & 2-\ve \\ \tfrac{1}{2}\ve & 2-\ve & 1+2\ve \\ 2-\ve & 1+2\ve & 6-6\ve \end{array}\right|=3+3\ve, 
$$
$$
 \Delta_4= \left|\begin{array}{cccc} 1+\ve & \tfrac{1}{2}\ve & 2-\ve & 1+2\ve \\ \tfrac{1}{2}\ve & 2-\ve & 1+2\ve & 6-6\ve \\ 2-\ve & 1+2\ve & 6-6\ve & 7+2\ve \\ 1+2\ve& 6-6\ve & 7+2\ve & 24-28\ve \end{array}\right|=7+10\ve, \ldots, 
$$
which (after shifting the index) gives the correct values of $X_1,X_2,X_3,X_4,X_5,\ldots$ for 
the combination of the original Somos-4 sequence (\ref{s4seq}) together with its odd part, 
namely  the shadow sequence $y_n^{(iv)}$ 
as in the fourth row of Table \ref{factortable}. 
\end{exa}

\begin{exa} For the recurrence (\ref{s4be}), the sequence with initial data given by four 1s with zero odd part was 
considered in \cite{ot}, which corresponds to $X_{-1}=X_0=X_1=X_2=1$ and  $\al^{(0)}=\be^{(0)}=\be^{(1)}=1$, 
$\al^{(1)}=0$ with our notation and indexing conventions. 
In this case, we take 
$$ 
\ha= 1+\ve, \quad \hb=1, \quad \hc=1+\tfrac{1}{2}\ve, \quad s_0=1, \quad s_1=-\tfrac{1}{2} \ve, 
$$
giving $U=-1$, $F=-3-\ve$, so that  
$$ 
\al=1, \qquad \be=1+\ve, \qquad J=4+\ve. 
$$  
Then the recursion (\ref{srec}) yields 
the sequence $(s_j)$ as 
$$ 
1, -\tfrac{1}{2} \ve, 2+\ve, 1-\ve, 6+4\ve,7,24+18\ve,41+18\ve, 115+98\ve, 
236+\tfrac{345}{2}\ve, \ldots , 
$$  
and 
the corresponding sequence of Hankel determinants 
begins with $\Delta_0=1$, $\Delta_1=1$, 
$$
\Delta_2
=\left|\begin{array}{cc} 1 & -\tfrac{1}{2}\ve \\ -\tfrac{1}{2}\ve & 2+\ve \end{array}\right|=2+\ve, 
\quad
\Delta_3= \left|\begin{array}{ccc} 1& -\tfrac{1}{2}\ve & 2+\ve \\ -\tfrac{1}{2}\ve & 2+\ve & 1-\ve \\ 2+\ve & 1-\ve & 6+4\ve \end{array}\right|=3+2\ve, 
$$
$$
 \Delta_4= \left|\begin{array}{cccc} 1 & -\tfrac{1}{2}\ve & 2+\ve & 1-\ve \\ -\tfrac{1}{2}\ve & 2+\ve & 1-\ve & 6+4\ve \\ 2+\ve & 1-\ve & 6+4\ve & 7 \\ 1-\ve& 6+4\ve & 7 & 24+18\ve \end{array}\right|=7+10\ve, \ldots, 
$$
of which the odd parts for index $n\geq 2$ give the (conjecturally) positive sequence of 
integers $1,2,10,48,160,1273,7346,51394,645078,...$ as found in \cite{ot}.
\end{exa} 

\begin{exa}
The analogous sequence for the recurrence (\ref{s4al}) was also presented in  
\cite{ot}, corresponding to $X_{-1}=X_0=X_1=X_2=1$ and  $\al^{(0)}=\be^{(0)}=\al^{(1)}=1$, 
$\be^{(1)}=0$. 
For this example, we take 
$$ 
\ha= 1+\ve, \quad \hb=1, \quad \hc=1+\tfrac{1}{2}\ve, \quad s_0=1, \quad s_1=0, 
$$
giving $U=-1-\tfrac{1}{2}\ve$, $F=-3-\ve$, so that  
$$ 
\al=1+\ve, \qquad \be=1, \qquad J=4+2\ve. 
$$  
Thus from the recursion (\ref{srec}) 
the sequence of dual number moments $(s_j)$ is found to be  
$$ 
1, 0, 2+\ve, 1+\tfrac{1}{2}\ve, 6+5\ve,7+\tfrac{13}{2}\ve,24+27\ve,41+\tfrac{105}{2}\ve, 115+164\ve, 
236+378\ve, \ldots , 
$$  
hence  
the corresponding sequence of Hankel determinants 
begins with $\Delta_0=1$, $\Delta_1=1$, 
$$
\Delta_2
=\left|\begin{array}{cc} 1 & 0 \\ 0  & 2+\ve \end{array}\right|=2+\ve, 
\quad
\Delta_3= \left|\begin{array}{ccc} 1& 0  & 2+\ve \\ 0 & 2+\ve & 1+\tfrac{1}{2}\ve \\ 2+\ve & 1+\tfrac{1}{2}\ve & 6+5\ve \end{array}\right|=3+3\ve, 
$$
$$
 \Delta_4= \left|\begin{array}{cccc} 1 &0  & 2+\ve & 1+\tfrac{1}{2}\ve \\ 0  & 2+\ve & 1+\tfrac{1}{2}\ve & 6+5\ve \\ 2+\ve & 1+\tfrac{1}{2}\ve & 6+5\ve & 7 +\tfrac{13}{2}\ve\\ 1+\tfrac{1}{2}\ve& 6+5\ve & 7 +\tfrac{13}{2}\ve& 24+27\ve \end{array}\right|=7+10\ve, \ldots, 
$$
so the integer sequence beginning with 
 $1,3,10,59,198,1387,9389,57983,752301,...$, as in \cite{ot}, 
is obtained from  the odd parts for index $n\geq 2$,  
and this 
is also conjectured to be positive.
\end{exa}

There is one more use for the $s_j$, obtained from the  moment generating function $G$ in (\ref{gcla}), 
that is relevant to shadow sequences, but now it is the classical 
case of $\C$-valued moments that concerns us. In that setting, the quantities $v_n$ appearing in the continued fraction (\ref{gcla}) can be written in the form 
\beq\label{vform}
v_n = \frac{\Delta_{n-1}^*}{ \Delta_{n-1}}
- \frac{\Delta_{n}^*}{ \Delta_{n}} \quad \mathrm{for}\,\,n\geq 1, 
\eeq 
where $\Delta_0^*=0$ and 
\beq\label{cas}
 \Delta_n^* = 
\left| \begin{array}{ccccc} 
s_0     & s_1 & \cdots & s_{n-2} & s_{n}    \\ 
s_1     &         & \iddots &    \vdots        & \vdots       \\
\vdots & \iddots &          &  \vdots       &   \vdots    \\ 
s_{n-2}     &  \cdots & \cdots & s_{2n-4} &  s_{2n-2}   \\
s_{n-1}     &  \cdots & \cdots & s_{2n-3} &  s_{2n-1}  
\end{array} 
\right|
\eeq  
is a bordered Hankel determinant for $n\geq 1$. 
Thus, as pointed out   in \cite{honecf}, the solution of the system (\ref{dtoda}) is expressed in terms of ratios 
of these Hankel/bordered Hankel determinants, whenever the moment sequence $(s_j)$ is generated 
by a recursion of the form (\ref{srec}) over $\C$. This yields a bordered Hankel  formula for the third shadow 
Somos-4 sequence.

\begin{propn} Let $(s_j)_{j\geq 0}$ be a moment sequence satisfying a recursion (\ref{srec}) 
over $\C$, corresponding to a particular 
Somos-4 sequence $(x_n)$, so that (up to rescaling)  it is given in terms of the 
associated Hankel determinants by  $x_n=\Delta_{n-1}$ for $n\geq 1$. Then, up to 
removing multiples of $y_n^{(i)}=x_n$, 
the third independent shadow Somos-4 sequence is given in terms of the corresponding bordered 
Hankel determinants (\ref{cas}) by   $y_n^{(iii)}=\Delta_{n-1}^*$  for $n\geq 1$.
\end{propn}  

\begin{prf}
Substitution of (\ref{vform}) into the formula (\ref{vsum}) gives a telescopic sum, 
which simplifies to 
$$ 
 y_n^{(iii)}=-x_n\left(v_0 -\frac{\Delta_{n-1}^*}{\Delta_{n-1}}\right),  
$$
and then using $x_n=\Delta_{n-1}$ this gives the required result after removing the multiple $-v_0 x_n$ 
from the front.
\end{prf}

\section{Conclusions} 
\setcounter{equation}{0}

We have given a very detailed decription of the solution of the initial value problem for 
the general dual Somos-4 recurrence (\ref{ds4}), from several different perspectives. The 
fact that this is possible relies heavily on the commutativity of the algebra of dual numbers, 
but  is also a reflection of the fact that there is an underlying discrete integrable system: 
either the QRT map (\ref{qrt}), which preserves the log-canonical symplectic form 
$\rd \log d_{n-1}\wedge \rd \log d_n$, or the map (\ref{dtoda}). (The latter map  
preserves a different symplectic structure, but the orbits of the two maps can be identified 
via a correspondence between 
their parameters and initial data.)
We have not explicitly described the analogues of these maps and their solutions over $\D$, but such 
a description is a straightforward consequence of our results on the dual Somos-4.

Continuous integrable systems with Grassmann variables have been studied for several decades, 
with one of the most recent results being a symmetry classification of $N=1$ supersymmetric 
scalar homogeneous evolutionary PDEs \cite{tw}. 
Many of our considerations here extend naturally to the dual number analogues of other 
discrete integrable systems, such as the family of maps in \cite{hkq}, which are connected with 
Gale-Robinson sequences and cluster algebras, or the higher genus analogues of (\ref{dtoda}) in \cite{honecf}. 
There are already versions of Yang-Baxter maps that include Grassmann variables, 
together with associated integrable lattice equations that satisfy the multidimensional 
consistency property \cite{kk}. 
Hopefully some of these techniques could also shed more light on superfriezes and cluster superalgebras, as in 
\cite{mgot, os}, which are relevant to Ptolemy relations in  super-Teichm\"uller theory \cite{moz}.

Finally, it would potentially be interesting to develop the  analytic theory of elliptic 
functions over the dual numbers, since the ubiquity of the derivatives $\partial_{g_2}$, $\partial_{g_3}$ 
suggests that modular identities might appear very naturally in this setting.

\noindent \textbf{Acknowledgments:} This research was supported by 
the  grant 
 IEC\textbackslash R3\textbackslash 193024 from the Royal Society. 

\end{document}